\documentclass{llncs}
\usepackage{amssymb}
\usepackage{amsfonts}
\usepackage{amssymb}
\usepackage{latexsym}
\usepackage{graphics}
\def\Hide#1{\relax}
\def\Fun#1#2{#1\rightarrow#2}
\def\Func#1{{\mathsf{#1}}}

\def\Recip#1{{#1}^{-1}}

\def\Max#1{\Func{max}{#1}}
\def\Min#1{\Func{min}{#1}}

\def\Abs#1{|#1|}
\def\Hom#1#2{\Func{Hom}(#1, #2)}
\def\Floor#1{\lfloor#1\rfloor}
\def\Ab{\Func{Ab}}
\def\Hom#1#2#3{\Func{Hom}_{#1}(#2, #3)}
\def\HomAb#1#2{\Hom{\Ab}{#1}{#2}}
\def\Sup{\Func{sup}}

\def\N{\mathbb{N}}

\def\E{\mathbb{E}}
\def\R{\mathbb{R}}
\def\Z{\mathbb{Z}}

\def\Rii{\sqrt{2}}
\def\Eu#1{{#1}^{\E}}

\newtheorem{Theorem}{Theorem}
\newtheorem{Lemma}[Theorem]{Lemma}

\newtheorem{Definition}[Theorem]{Definition}
\def\Proof{\par {\bf Proof: }}
\def\Done{\rule{0.75em}{0.75em}}
\title{The Eudoxus Real Numbers}
\author{R.D. Arthan}
\institute{Lemma 1 Ltd. \\ 2nd Floor, 31A Chain Street, \\ Reading UK \hspace{1em} RG1 2HX\\
\email{rda@lemma-one.com}}
\begin{document}
\maketitle
\begin{abstract}
This note describes a representation of the real numbers due to Schanuel.
The representation lets us construct the real numbers from first principles.
Like the well-known construction of the real numbers using Dedekind cuts, the idea is inspired
by the ancient Greek theory of proportion, due to Eudoxus.
However, unlike the Dedekind construction, the construction proceeds directly from the
integers to the real numbers bypassing the intermediate construction of the rational numbers.

The construction of the additive group of the reals depends on rather
simple algebraic properties of the integers. This part of the construction can be generalised
in several natural ways, e.g., with an arbitrary abelian group playing the role of the integers.
This raises the question: what does the construction construct?
In an appendix we sketch some generalisations and answer this question in some simple cases.

The treatment of the main construction
is intended to be self-contained and assumes familiarity only
with elementary algebra in the ring of integers and with the definitions of the abstract algebraic notions of group, ring and field.

\end{abstract}
\section{Introducing $\E$}

\begin{quotation}
{\it
Magnitudes are said to be {\em in the same ratio}, the first to the second and the third to the fourth, when, if any equimultiples whatever are taken of the first and third, and any equimultiples whatever of the second and fourth, the former equimultiples alike exceed, are alike equal to, or alike fall short of, the latter equimultiples respectively taken in corresponding order. \\
Euclid. Elements of Geometry. Book V. Definition 5.
}
\end{quotation}

That is to say the ratios $a:b$ and $c:d$ are equal iff. for all positive integers $m$ and $n$,
the statements $ma > nb$ and $mc > nd$ are either both true or both false, and similarly
for $ma = nb$ and $ma < nb$.

\begin{figure}
\begin{center}
\includegraphics{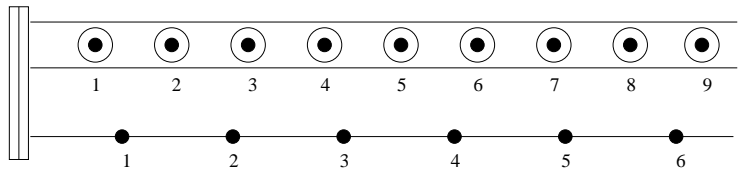}
\caption{De Morgan's Colonnade and Fence}
\label{Colonnade}.
\end{center}
\end{figure}

According to the commentary in Heath's translation of the Elements \cite{Heath56},
de Morgan gave an interesting rationale for this definition:
imagine a fence with equally spaced railings in front of a colonnade of equally spaced columns,
as shown, in plan, in figure~\ref{Colonnade}. Let the distance between the columns be $C$ and the distance between
the railings be $R$. If the construction is continued indefinitely, 
an observer can compare $C$ with $R$ to any degree of accuracy without making any measurements
just by counting the columns and railings.
For example, in the figure, the 6th railing lies
between the 8th and 9th columns, so that $8C < 6R < 9C$ which means that $R/C$ lies
between $4/3$ and $3/2$.
If more precision were required, the observer might continue counting to find that the 25th railing lies between the 35th and 36th columns and conclude that $R/C$ lies between $7/5$ and $36/25$.

This picture suggests a way of representing real numbers:
construct the colonnade so that  the distance $C$ between the columns is $1$
and represent $R$ by the sequence of integers
in which the $m$-th term, $R_m$ say, gives the number of columns
to the left of or in line with the $m$-th railing in the figure. This sequence $R_m$ will represent $R$.
In the example in figure~\ref{Colonnade}, the first few $R_m$ are  $1, 2, 4, 5, 7, 8, \ldots$.
Since Euclid's book V is generally believed to describe the work of Eudoxus, let us call the real numbers
represented in this way the {\em Eudoxus reals}.

We can allow a bounded amount of error in the calculation of the $R_m$ without affecting
the number represented. For example, if we consistently forget to count the first column, the
resulting sequence $R'_m = R_m - 1$ should still represent $R$, since for large $m$ the difference
between $R'_m/m$ and $R_m/m$ will be negligible.

If $R$ happens to be a positive integer, then the sequence $R_m = mR$ will be an additive homomorphism,
i.e., the equation $R_{m+n} = R_m + R_n$ will hold for all $m$ and $n$. For arbitrary $R$, this equation
should be ``approximately true''.
The observations of the previous paragraph suggest that the right notion of ``approximately true''
is true with a bounded amount of error. This leads to the following definition:

\begin{Definition}\label{AH}
A function $f$ from $\Z$ to $\Z$ is said to be an {\em almost homomorphism} iff. the function
$d_f$ from $\Z \times \Z$ to $\Z$ defined by  $d_f(p, q) =f(p+q) - f(p) - f(q)$
 has bounded range, i.e., for some integer $C$, $\Abs{d_f(p, q)}.< C$ for all $p, q \in \Z$.
\end{Definition}

The function $d_f$ measures the extent to which $f$ fails to be a genuine homomorphism:
$f$ is a genuine homomorphism iff. $d_f$ is identically $0$.
Note that the requirement on an almost homomorphism is {\em not} that $d_f(p, q)$ be
non-zero for only finitely many $p$ and $q$, but that values taken on by $d_f(p, q)$ 
be bounded as $p$ and $q$ range over $\Z$.  Typically, if $d_f(p, q)$ is non-zero
for some $p$ and $q$, there will be infinitely many pairs $(r, s)$ for which $d_f(r, s) = d_f(p, q)$.

The set $\Fun{\Z}{\Z}$ of all functions from $\Z$ to $\Z$ becomes an abelian group if we
add and invert functions pointwise: $(f+g)(p) = f(p) + g(p)$, $(-f)(p) = -f(p)$.
It is easily checked that if $f$ and $g$ are almost homomorphisms then so are $f+g$ and $-f$.
Writing $S$ for the set of all almost homomorphisms from $\Z$ to $\Z$, this
shows that $S$ is a subgroup of $\Fun{\Z}{\Z}$.
Let us write $B$ for the set of all functions from $\Z$ to $\Z$ whose range is bounded.
$B$ is a subgroup of $S$. 

\begin{Definition}\label{Edef}
The abelian group $\E$ of {\em Eudoxus reals} is the quotient group $S/B$.
\end{Definition}

Let me spell out what this means: elements of $\E$ are equivalence classes, $[f]$ say, where $f$ is
an almost homomorphism from $\Z$ to $\Z$, i.e., $f$ is a function from $\Z$ to $\Z$ such
that $d_f(p, q) = f(p, q) -f(p) - f(q)$ defines a function from $\Z \times \Z$ to $\Z$ whose
range is bounded. We have $[f] = [g]$ iff.
the difference $f - g$ has bounded range, i.e., iff. $\Abs{f(p) - g(p)} < C$ for some $C$ and all $p$ in $\Z$.
The addition and inverse in $\E$ are induced by the pointwise addition and inverse of representing almost
homomorphisms: $[f] + [g] = [f + g]$, $-[f] = [-f]$ where $f+g$ and $-f$ are defined by $(f+g)(p) = f(p) + g(p)$ and $(-f)(p) = -f(p)$ for
all $p$ in $\Z$.

For example, if integers $A$ and $B$ are given, the linear function $f(p) = Ap + B$ is an almost
homomorphism, with $d_f(p, q) = - B$.
If $g(p) = Cp + D$ is another such function then $[f] = [g]$ iff. $A = C$.
It follows that we can define a one-one mapping $\Eu{\_}$ from $\Z$ into $\E$ by setting
$\Eu{A} = [p \mapsto Ap]$. This mapping is a homomorphism: $\Eu{(A+B)} = \Eu{A}+\Eu{B}$.

Our aim is now to show that the group $\E$ of Eudoxus reals is isomorphic to the group of real numbers
under addition. To do this, we will prove
``from first principles'', i.e., using only elementary set theory and integer arithmetic, that, like $\R$, $\E$ can be
equipped with an ordering and a multiplication that make it into a complete ordered field. Any two complete ordered fields are isomorphic and
so $\E$ and $\R$ (howsoever $\R$ be constructed) are isomorphic.

\section{Ordering}
We will start by investigating the ordering of $\E$. In doing so we prove some facts
about the asymptotic behaviour of almost homomorphisms which will be useful in the
sequel. 

Our first lemma gives a partial lower bound on the growth rate of certain almost homomorphisms.
In terms of de Morgan's colonnade and fence, it says that by picking a large number $M$
and removing all the railings except the ones numbered $M$, $2M$, $3M$ and so on, the rate
of growth of the $R_m$ can be made to exceed any given linear bound.

In the proof of this lemma, we use the identity  $f(p+q) = f(p) + f(q) + d_f(p, q)$ satisfied
by the function $d_f$ of definition~\ref{AH}. We will use this identity and its elementary consequences
without further mention in the sequel.

\begin{Lemma}\label{LowerBoundLemma}
If $f$ is an almost homomorphism such that $f(m)$ takes on infinitely many positive
values as $m$ ranges over $\N$, then, for any $D > 0$, there
is an $M > 0$ such that $f(mM) > (m+1)D$ for all positive
integers $m$.
\end{Lemma}
\Proof
By the definition of almost homomorphism, there is a $C$ such that
for any $p$ and $q$ in $\Z$, $\Abs{d_f(p, q)} < C$.
Put $E = C + D$.
Since $f(m)$ takes on infinitely many positive values as $m$ ranges over $\N$, we may
choose $M > 0$ such that $f(M) > 2E$.
Assume, by induction, that $f(mM) > (m+1)E$ for some positive integer $m$.
We then have $f((m+1)M) = f(mM) + f(M) + d_f(mM, M) > f(mM) + f(M) - E > (m+1)E + 2E - E = (m+2)E$.
Thus, by induction $f(mM) > (m+1)E$ for all $m > 0$.
But then $f(mM) > (m+1)E = (m+1)(C + D) > (m+1)D$, and the lemma holds.
\Done

\begin{Definition}\label{OrderedGroupDef}
An {\em ordered abelian group} is an abelian group $G$, together with a set $P \subseteq G$ of {\em positive elements} satisfying the following:
\begin{enumerate}
\item the identity, $0$, of the group is not a member of $P$;
\item $P$ is closed under addition, i.e., if $x \in P$ and $y \in P$, then $x+y \in P$;
\item If $x$ is any non-zero element of $G$ then either $x$ or $-x$ is a member of $P$.
\end{enumerate}
\end{Definition}

If $G$ is an ordered abelian group with positive elements $P$, and $x$ and $y$ are elements of $G$
we write $x < y$ for $y - x \in P$ (and use $x \le y$, $x \ge y$ and $x > y$ in the usual way,
to mean $\lnot y < x$, $y \le x$ and $y < x$ respectively).
One may check that $<$ is a total ordering
on the elements of $G$ which is compatible with addition
in the sense that whenever $x < y$ then also $x + z < y + z$ for any $z$.
The absolute value $\Abs{x}$ of an element $x$ of $G$ is defined to be the larger of $x$ and $-x$
(or $0$ if $x = 0$). Absolute values satisfy the triangle inequality: $\Abs{x+y} \le \Abs{x} + \Abs{y}$ for any $x$ and
$y$ in $G$.

The following important lemma will justify our definition of the ordering of the Eudoxus reals.
For $x = [f] \in \E$, the three cases of the lemma will correspond to $x = 0$, $x > 0$ and $x < 0$ respectively.

\begin{Lemma}\label{TrichLemma}
Let $f$ be an almost homomorphism, then exactly one of the following holds:
\begin{enumerate}
\item $f$ has bounded range;
\item for any $C > 0$, there is an $N$, such that $f(p) > C$ whenever $p > N$;
\item for any $C > 0$, there is an $N$, such that $f(p) < -C$ whenever $p > N$.
\end{enumerate}
\end{Lemma}
\Proof
The three cases are clearly mutually exclusive. To prove that at least one of the cases holds
for any almost homomorphism $f$, let us
first assume that $f(m)$ takes on infinitely many positive values for $m \in \N$.
Our plan is to relate $f$ to a function $g$ whose growth rate is more tractable.
Lemma~\ref{LowerBoundLemma} gives lower bounds on the growth of the
values $f(mM)$ for certain $M$.  $g(p)$ will agree with $f(p)$ when $p$  has the form $mM$ and will
``mark time'' between values of this form.

To define $g$,  first choose
 $D$ so that $\Abs{d_f(p, q)} < D$ for all $p$ and $q$. 
By lemma~\ref{LowerBoundLemma} there is a positive integer $M$
such that $f(mM) > (m+1)D$ for all positive integers $m$.
Given any integer $p$, there are unique integers $d$ and $r$ such that $p = dM + r$ and $0\le r < M$,
and given such $d$ and $r$ we define $g(p) = f(dM)$.

Choose $E$ such that $E > \Abs{f(r)}$ for $0 \le r < M$.
Given $p \in \Z$, let $d, r \in \Z$ be such that $p = dM + r$ and $0 \le r < M$. We then have $\Abs{(f - g)(p)} = \Abs{(f - g)(dM+r)} = \Abs{f(dM + r) - g (dM + r)} = \Abs{f(dM) + f(r) + d_f(dM, r) - f(dM)} = \Abs{f(r) + d_f(dM, r)} < E + D$.
So $\Abs{(f - g)(p)}$ is bounded for $p \in \Z$ with bound $B = E + D$.

Given $C$, pick $n > 0$ such that $(n+1)D > B+C$ and set $N = nM$.
If $p >  N = nM$, then, when we write $p = dM + r$ with $0 \le r < M$, we
must have $d \ge n$, but then $g(p) = f(dM) > (d+1)D \ge (n+1)D > B+C$, so that $f(p) > g(p) - B >  C$.
Thus case {\em(ii)} holds if $f(m)$ takes on infinitely many positive values for $m \in \N$.

If $f(m)$ takes on infinitely many negative values for $m \in \N$,
then by the above argument, case {\em(ii)} holds for $-f$, but then case {\em(iii)} holds for $f$.

Finally, if $f(m)$ is bounded for $m \in \N$, then $f$ is bounded, because for $p < 0$,
$f(p) = f(0) - f(-p) - d_f(-p, p)$, so that, as $d_f(-p, p)$ is bounded, negative values
of $p$ contribute only a finite number of values to the range of $f$.
Thus case {\em(i)} applies.
\Done

We now define the set $P$ of positive elements of the group $\E$
to comprise all elements $[f]$ such that for
$m \in \N$, $f(m)$ takes on infinitely many positive values, so that case {\em(ii)}
of the above lemma applies to $f$.
Using lemma~\ref{TrichLemma}, it is easy to see that this definition is independent
of the choice of representative $f$. Again by lemma~\ref{TrichLemma}, the set $P$ does not contain $0$,
and is such that for any $[f]$ in $\E$, either $[f] \in P$ or $-[f] \in P$.

\begin{Theorem}\label{OrderedGroupThm}
The group, $\E$, of Eudoxus reals becomes an ordered abelian group if we take the set, $P$, of positive
elements to comprise all elements $[f]$ for which $f(m)$ takes on infinitely many positive
values as $m$ ranges over $\N$.
\end{Theorem}
\Proof
From the remarks above, all we have to show is that the set $P$ of positive elements is closed under
addition. So, let $x = [f]$ and $y = [g]$ be members of $P$, so that both $f$ and $g$ fall
under case {\em(ii)} of lemma~\ref{TrichLemma}, i.e., given $C > 0$, there are integers $M$ and $N$ such that
$f(m) > C$ whenever $m > M$ and $g(m) > C$ whenever $m > N$;
but then, if $m > \Max{\{M, N\}}$,  we will have $(f+g)(m) = f(m) + g(m) >  C + C > C$,
so that $(f+g)(m)$ takes on infinitely many positive values as $m$ ranges over $\N$ 
and $x+y = [f+g]$ belongs to $P$ as required.
\Done

\section{Multiplication}
In terms of de Morgan's colonnade and fence analogy,
if we are given fences representing numbers $R$ and $S$ say, the fence representing the product $RS$ can
be constructed by scaling the fence for $R$ to give a new fence whose railings bear the same
relation to the railings in the fence for $S$ as those in the fence for $R$ do to the columns in the colonnade.
If $T_m$ is the number of columns corresponding to the $m$-th railing in this new fence,
we would have $T_m = S_{R_m}$. Of course, we expect multiplication to be commutative, i.e.,
we would hope that $S_{R_m} = R_{S_m}$, but this is not true in general. However, it turns out
to be ``almost true'', and that will be good enough for our purposes.

The above remarks suggest that  multiplication should correspond to composition of
almost homomorphisms. 
So let us try to make $\E$ into a ring by
defining the product $[f][g]$ to be  $[f \circ g] = [p \mapsto f(g(p))]$.
It is easy to check that $[f][g]$ is well-defined, i.e., that $f \circ g$ is an almost homomorphism if $f$ and $g$ are, and that $[f \circ g]$ is independent of the choice of representatives $f$ and $g$.
It is also easy to see that the identity function $1 = \Eu{1}$  satisfies $1x = x1 = x$
and that $x(yz) = (xy)z$ and $(x +y)z = xy + xz$ for all $x$, $y$ and $z$ in $\E$.
It is not so easy to see that $x(y + z) = xy + xz$, but this will follow, and so $\E$ will be a ring,
if we can show that the multiplication is commutative. We will do this using the following lemma.

In the proof of the lemma, and in later work, we talk of ``adding'' two inequalities
of the form $\Abs{P} < A$ and $\Abs{Q} < B$, by which we mean adding the
inequalities in the usual sense to give $\Abs{P}+\Abs{Q} < A + B$ and then using the triangle inequality
$\Abs{P+Q} \le \Abs{P} + \Abs{Q}$ to conclude $\Abs{P+Q} < A + B$.

\begin{Lemma}\label{MultLemma}
If $f$ is an almost homomorphism, say with $\Abs{d_f(p, q)} < C$ for all $p$ and $q$,
then the following holds for all $p$ and $q$ in $\Z$:
\[
\Abs{pf(q) - qf(p)} < (\Abs{p} + \Abs{q} + 2)C
\]
\end{Lemma}
\Proof
I claim that $\Abs{f(pq) - pf(q)} < (\Abs{p} + 1)C$ for all $p$ and $q$.
We have $\Abs{f(0)} = \Abs{f(0+0) - f(0) - f(0)} = \Abs{d_f(0, 0)} < C$, so the claim
is certainly true when $p = 0$. For $p \ge 0$, 
we have $\Abs{f((p+1)q) - f(pq)  - f(q)} = \Abs{d_f(pq, q)}< C$;
if  we have $\Abs{f(pq) - pf(q)} < (\Abs{p} + 1)C$,  then adding these inequalities
gives $\Abs{f((p+1)q) - (p+1)f(q)} < (\Abs{p} + 2)C$, so the claim is true for all $p \ge 0$
by induction. Similarly, using $\Abs{-f(pq) + f(-(p+1)q) + f(q)} = \Abs{-d_f(-(p+1)q, q)}$,
we have that the claim is true for all $p < 0$ too.

We now know that $\Abs{f(pq) - pf(q)} \le (\Abs{p} + 1)C$ for all $p$ and $q$.
Interchanging $p$ and $q$ and rearranging gives that $\Abs{qf(p) - f(pq)} < (\Abs{q} + 1)C$
Adding these inequalities gives $\Abs{pf(q) - qf(p)} < (\Abs{p} + \Abs{q} + 2)C$ as required.
\Done.

To introduce the ordering on $\E$ we used linear lower bounds on the growth rates of
almost homomorphisms. To show that multiplication is commutative, we will
need the linear upper bounds supplied by the following lemma.

\begin{Lemma}\label{UpperBoundLemma}
If $f$ is an almost homomorphism, there exist integers $A$ and $B$ such that for any $p$ in $\Z$
we have $\Abs{f(p)} < A\Abs{p} + B$;
\end{Lemma}
\Proof
Let $C$ be such that $\Abs{d_f(p, q)} < C$ for all $p$ and $q$.
Taking $q =1$ in lemma~\ref{MultLemma} gives $\Abs{f(p) - pf(1)} < (\Abs{p} + 3)C$,
so that $\Abs{f(p)} < (\Abs{p} + 3)C + \Abs{p}\Abs{f(1)} = (C+\Abs{f(1)})\Abs{p} + 3C$ and we may take $A = C + \Abs{f(1)}$
and $B =3C$.
\Done

\begin{Theorem}
The multiplication on $\E$ induced by composition of almost homomorphisms makes
$\E$ into a commutative ring with unit.
\end{Theorem}
\Proof
By the remarks preceding lemma~\ref{MultLemma} above, all we have left to do is to show that the multiplication is commutative.
This means that, given almost homomorphisms $f$ and $g$, we must show that the function $(f o g) - (g o f)$
has bounded range. To see this, take $q = g(p)$ in the inequality of lemma~\ref{MultLemma}.
This gives:
\[
	\Abs{pf(g(p)) - g(p)f(p)} < (\Abs{p} + \Abs{g(p)} + 2)C
\]
and, similarly, interchanging $f$ and $g$ and rearranging, we have:
\[
	\Abs{g(p)f(p) -p(g(f(p))} < (\Abs{p} + \Abs{f(p)} + 2)C
\]
Adding these two inequalities gives:
\[
	\Abs{pf(g(p)) - p (g(f(p))} < (2\Abs{p} + \Abs{f(p)} + \Abs{g(p)} + 2)C
\]
Lemma~\ref{UpperBoundLemma} gives us linear upper bounds for $\Abs{f(p)}$ and $\Abs{g(p)}$.
Using these bounds, we can find integers D and E such that
\[
	\Abs{p} \Abs{f(g(p)) - (g(f(p))} < D \Abs{p} + E
\]
But then for all sufficiently large values of $\Abs{p}$, we must have, 
\[
	\Abs{f(g(p)) - (g(f(p))} < D + 1
\]
So that $(f o g) - (g o f)$ does indeed have bounded range, multiplication is commutative
in $\E$ and $\E$ is a commutative ring with unit $\Eu{1}$.
\Done

The following simple lemma eases the burden of showing that
some functions are almost homomorphisms.

\begin{Lemma}\label{ExtensionLemma}
If $f$ is a function from $\Z$ to $\Z$ such that, {\em(i)}, $f(p) = -f(-p)$ whenever $p < 0$,
and, {\em(ii)}, $d_f(m, n) = f(m+n) - f(m) - f(n)$ is bounded for $m$ and $n$ in $\N$,
then $f$ is an almost homomorphism.
\end{Lemma}
\Proof
We have to show that $d_f(p, q)$ takes
on only finitely many values as $p$ and $q$ range over $\Z$.
We are given that $d_f(p, q)$ takes on only finitely many values if $p \ge 0$ and $q \ge 0$.
If $p$ and $q$ are both negative, then
$d_f(p, q) = -d_f(-p, -q)$, so that this case contributes only finitely many extra values for $d_f(p, q)$.
Finally, assume exactly one of $p$ or $q$ is negative, say $p$: if $p + q \le  0$, put $a = q$
and $b = -(p+q)$;
if $p + q >0$ put $a = p + q$ and $b = -p$; in both cases we have that $a \ge 0 $ and $b \ge 0$
and we find $\Abs{d_f(p, q)} = \Abs{d_f(a, b)}$ so that the case when exactly one of $p$ or $q$ is negative contributes no additional values and $f$ is indeed an almost homomorphism.
\Done

We now show that $\E$ is actually a field: this means that any non-zero element  has a multiplicative
inverse. In terms of de Morgan's colonnade and fence, to construct
the inverse, we would just interchange the roles of the columns and the railings. The construction
in the proof below formalises this.

\begin{Theorem}
For any non-zero $x$ in $\E$, there is an element $\Recip{x}$ such that $x\Recip{x} = \Eu{1}$.
Thus the commutative ring $\E$ is a field.
\end{Theorem}
\Proof
First let us assume that $x = [f]$ is a positive element of $\E$. By lemma~\ref{TrichLemma},
for all $m \in \N$, $f(n) > m$ for all sufficiently large $n$, and so we may
define a function $g$ from $\Z$ to $\Z$ as follows:
\[
g(p) =
 \left\{
    \begin{array}{ll}
	\Min \{n : \N | f(n) \ge p\} & \mbox{if $p \ge 0$} \\
	-g(-p)	& \mbox{otherwise}
    \end{array}
 \right.
\]

I claim $g$ is an almost homomorphism. By lemma~\ref{ExtensionLemma}, it suffices
to check that $d_g(m, n)$ is bounded for $m, n \in \N$.
To see this note that for all but finitely many $m, n \in \N$, $p = g(m)$ and
$q = g(n)$ are both positive, and then certainly $r = g(m+n)$ is also positive.
What we have to prove is that $d_g(m, n) = r - p - q$ is bounded as $m$ and $n$ range over $\N$.
By the definition of $g$, if $m$ and $n$ are large enough for $p$ and $q$ to be positive, we have:
\[
\begin{array}{rcccl}
	f(p) &\ge& m &>& f (p-1) \\
	f(q) &\ge& n &>&  f(q-1) \\
	f(r) &\ge& m+n &>& f(r-1)
\end{array}
\]

From these inequalities, we can derive:
\[\begin{array}{rcccl}
	f(r) - f(p-1) - f(q-1) &>& (m+n) - m - n &=& 0 \\
	f(r-1) - f(p) - f(q)  &<& (m+n) - m - n &=& 0 \\
\end{array}
\]
For each of the above two inequalities, the difference between the left-hand side and
$f(r - p - q)$ is bounded (independently of $p$, $q$ and $r$) because
$f$ is an almost homomorphism.
It follows that $f(r - p - q)$ is bounded as $m$ and $n$ range over $\N$,
but then $r - p - q$ must be be bounded, since we
are assuming that $f$ is positive, so that if $t$ ranges over an unbounded set of integers, so also does $f(t)$,
by lemma~\ref{TrichLemma} applied to $f$ and $(p \mapsto f(-p))$ (using the fact that $f(p) + f(-p)$ is bounded).  
Thus $d_g(m, n) = r - p - q$ is bounded as $m$ and $n$ range over $\N$,
and so, by lemma~\ref{ExtensionLemma}, $g$ is indeed an almost homomorphism.

For large enough $m$, we have
\[
	f(g(m)) \ge m > f(g(m) - 1) \ge f(g(m)) - C
\]
where $C$ is independent of $m$. Since $\Eu{1} = (m \mapsto m)$, and $g$ is
an almost homomorphism, it follows that
$[f][g] - \Eu{1} = [f \circ g] - \Eu{1} = \Eu{0}$. I.e., $[f][g] = \Eu{1}$.
Thus $[g]$ is a multiplicative inverse for the positive element $x = [f]$.

If $x = [f]$ is negative, apply the above construction to give an inverse, $[g]$ say, for the
positive element $-[f]$, we then have $[f](-[g]) = (-[f])[g] = \Eu{1}$, so $-[g]$ provides
the required inverse for $x = [f]$.
\Done

We note in passing that the almost homomorphism $g$ constructed above to represent
the inverse of a positive $x$ is non-decreasing:
i.e., for any $p$, $g(p) \le g(p+1)$. Thus if we repeat the construction to give a representative
for $\Recip{[g]}$, we get a non-decreasing almost homomorphism, $h$ say, such that $[h] = \Recip{({\Recip{x}})} = x$.
Thus any positive (resp. negative) element of $\E$ can be represented by a non-decreasing (resp. non-increasing) almost homomorphism.

\begin{Definition}
An {\em ordered field} is a field $F$ together with a set $P \subseteq F$ of {\em positive elements}
which makes the additive group of $F$  an ordered group as in definition~\ref{OrderedGroupDef}
and satisfies the following:
\begin{enumerate}
\item the unit element, $1$, is in $P$
\item $P$ is closed under multiplication, i.e., if $x \in P$ and $y \in P$, then $xy \in P$.
\end{enumerate}
\end{Definition}

If $F$ is an ordered field, we define $x < y$, $x \le y$ etc. just as we did for an ordered group.
The ordering is then compatible with the multiplication in the sense that multiplication by
positive elements is strictly order-preserving (if $x > 0$ and $y < z$, then $xy < xz$)
 and multiplication by negative elements is strictly order-reversing
(if $x < 0$ and $y < z$, then $xy > xz$).
A non-zero element $x$ of $F$ is positive iff. $\Recip{x}$ is positive
and taking multiplicative inverses of positive elements is strictly order-reversing
(if $0 < x < y$ then $\Recip{x} > \Recip{y} > 0$).

\begin{Theorem}\label{OrderedFieldThm}
Under the ordering of theorem~\ref{OrderedGroupThm},
the field $\E$ is an ordered field.
\end{Theorem}
\Proof
Evidently the unit element $1 = \Eu{1}$ is positive and so, given theorem~\ref{OrderedGroupThm}, 
we have only to show that the set of positive elements of $\E$ is closed under
multiplication. By lemma~\ref{TrichLemma}, an element $x = [f]$ of $E$ is positive
iff. for any $C$ there is an $M$ such that $f(p) > C$ whenever $p > M$.
Assume $x = [f]$ and $y = [g]$ are positive. Given $C$, choose $M$ such
that $f(p) > C$ whenever $p > M$ and then choose $N$ such that $g(p) > M$
whenever $p > N$. Then if $p >  N$, $f(g(p)) > C$ so that $xy = [f \circ g]$
is positive as required.
\Done

If $F$ is any ordered field, there is a unique order-preserving homomorphism, $h$ say, from $\Z$  to $F$
(defined inductively by $h(0) = 0_F$, $h(p+1) = h(p) + 1_F$ and $h(p-1) = h(p) - 1_F$, where
$0_F$ and $1_F$ are the zero and unit elements of $F$ respectively).
The function $p \mapsto \Eu{p}$, where $\Eu{p} = [q \mapsto pq]$ gives an explicit
representation of this order-preserving homomorphism for the ordered field $\E$. 
For $x, y \in \E$ we will write $x/y$ for $x\Recip{y}$ and $\Abs{x}$ for the absolute value of $x$, i.e., the
larger of $x$ and $-x$.  Clearly $\Abs{\Eu{p}} = \Eu{\Abs{p}}$.

We can now give algebraic interpretations of lemmas~\ref{LowerBoundLemma}
and~\ref{UpperBoundLemma}.
Lemma~\ref{LowerBoundLemma} says that for
any positive $x$ in $\E$ and for all positive $D$, there is
a positive integer $M$ such that $x\Eu{M} \ge \Eu{D}$, i.e., $x \ge \Eu{D}/\Eu{M}$.
(to see this, put $x = [f]$ and apply the lemma to $D$ giving an $M$ such that
$[f]\Eu{M} = [p\mapsto f(pM)] \ge [p \mapsto (p+1)D] = [p \mapsto pD] = \Eu{D}$).
Similarly lemma~\ref{UpperBoundLemma} says that for any $x$ in $\E$, there
is an integer $A$ such that $\Eu{A} \ge x$.
A slight strengthening and some consequences of this are given in the following lemma:

\begin{Lemma}\label{ArchLemma}
\begin{enumerate}
\item
If $x$ is any element of $\E$, then there is a positive integer $M$ such that
$x < \Eu{M}$.
\item
The ordered field $\E$ is archimedean, i.e.,
if  $x, y$ are elements of $\E$, with $y$ positive, then there is a positive integer $M$ such that
$x < \Eu{M}y$.
\item
If $x < y$ are any elements of $\E$, then there are integers $M$, $N$, with $N$ positive
such that $x < \Eu{M}/\Eu{N} < y$.
\item
The ordering on $\E$ is dense: i.e., if $x < y$, then there is a $z$ such that $x < z < y$.
\end{enumerate}
\end{Lemma}
\Proof
For part {\em(i)},
By the remarks above, there is an integer $A$ such that $x \le \Eu{A}$, so we may take $M = \Max{\{1, A + 1\}}$.

For part {\em(ii)},
by part {\em(i)}, there are positive integers $C$ and $D$
such that $x < \Eu{C}$ and $\Recip{y} < \Eu{D}$.
But then multiplying the latter through by the positive element $y$, we have $1 <\Eu{D}y$,
whence $x < \Eu{C} < \Eu{C}\Eu{D}y = \Eu{(CD)}y$, and we may take $M = CD$.

For part {\em(iii)}, by part {\em(i)},
there is a positive integer $N$ such that $1/(y-x) < \Eu{N}$.  But then
as $0 < y - x$, we have $0 < 1/\Eu{N} < y - x$, or, equivalently, $x < x + 1/\Eu{N} < y$.
Applying part {\em(ii)} to $\Abs{x}$ and $1/\Eu{N}$, there is a $K$ such that
$\Abs{x} < \Eu{K}/\Eu{N}$, whence
$-\Eu{K}/\Eu{N} <  x < \Eu{K}/\Eu{N}$.
Let $L$ be the largest integer such that $\Eu{L}/\Eu{N} \le x$ and set $M = L + 1$.
Then we have $x < \Eu{M}/\Eu{N}$ by the choice of $L$.
Also, we must have $\Eu{M}/\Eu{N} < y$, since if not, then $y \le \Eu{M}/\Eu{N}$  and we would have
$1/\Eu{N} < y - x \le \Eu{M}/\Eu{N} - x$ which gives $x < \Eu{(M-1)}/\Eu{N} = \Eu{L}/\Eu{N}$
contradicting our choice of $L$.

Part {\em(iv)} is an immediate consequence of part {\em(iii)}.
\Done

\section{Completeness}

If $x$ is any element of $\E$, then applying part {\em(i)} of lemma~\ref{ArchLemma} to $\Abs{x}$ gives an integer $M$ such
that $\Eu{(-M)} < x < \Eu{M}$. So the set of integers $p$ such that $\Eu{p} \le x$ is
non-empty and bounded above. This justifies the following definition:

\begin{Definition}\label{FloorDef}
If $x$ is any element of $\E$, define the {\em floor} of $x$, written $\Floor{x}$,
by the conditions that $\Floor{x} \in \Z$ and $\Eu{\Floor{x}} \le x < \Eu{(\Floor{x} + 1)}$.
\end{Definition}

The floor function from $\E$ to
$\Z$ is weakly order-preserving, i.e., if $x \le y$ in $\E$ then $\Floor{x} \le \Floor{y}$
in $\Z$,  and is a left inverse to $\Eu{p}$, i.e., $\Floor{\Eu{p}} = p$.
Floor does preserve some strict inequalities: if $x, y \in E$, $p \in \Z$ and $x < \Eu{p} \le y$,
then $\Floor{x} < \Floor{y}$ and, in particular, $\Floor{x} < \Floor{\Eu{p}} = p$.
We also have the identity $\Floor{\Eu{p} + x} = p + \Floor{x}$ for any $p \in \Z$ and $x \in \E$.

The floor function is not a homomorphism. However, it is an almost homomorphism, if we generalise
definition~\ref{AH} in the obvious way to cover functions from $\E$ to $\Z$:

\begin{Lemma}\label{FloorAHLemma}
$\Floor{x+y} - \Floor{x} - \Floor{x}$ ranges over a bounded subset of $\Z$, namely, $\{ 0, 1\}$, as $x$ and $y$ range over $\E$.
\end{Lemma}
\Proof
We have the following three inequalities in $\E$, where the third is obtained by adding the first two
and simplifying using the fact that the function $\Eu{\_}$ is a  homomorphism

\[
\begin{array}{rcl}
\Eu{0} &\le& x - \Eu{\Floor{x}}< \Eu{1} \\
\Eu{0} &\le& y - \Eu{\Floor{y}}< \Eu{1} \\
\Eu{0} &\le& x + y - \Eu{\Floor{x}} - \Eu{\Floor{y}} < \Eu{2}
\end{array}
\]
Applying the function $\Floor{\_}$ to the third inequality and using the remarks following
definition~\ref{FloorDef}, we then have the following inequality in $\Z$:
\[
\begin{array}{rcl}
0 &\le& \Floor{x+y} - \Floor{x} - \Floor{y} < 2
\end{array}
\]
and that completes the proof.
\Done

We will now show that the ordering on $\E$ is complete in the sense of the following definition.

\begin{Definition}
A totally ordered set $A$ is {\em complete} iff.
any $S \subseteq A$ that is non-empty and bounded above has a {\em supremum} in $A$,
i.e. an element $s$ of $A$ (not necessarily of $S$) such that
\begin{enumerate}
\item for any $x$ in $S$, $x \le s$
\item for any $y$ in $A$, if $x \le y$ whenever $x \in S$, then $s \le y$.
\end{enumerate}
\end{Definition}

We must show that every non-empty, bounded above subset $S$ of $\E$ has a supremum.
The construction may be understood in terms of de Morgan's colonnade and fence picture as follows:
we now have a set of fences;  we know that for each $m$, there is a column that is to the right
of the $m$-th railing in all of the fences; we construct a new fence in which the $m$-th railing
is in line with the right-most column that is to the left of or in line with the $m$-th railing in 
at least one of the fences. This new fence represents the supremum of the set of numbers
represented by the original fences.

\begin{Theorem}\label{CompleteThm}
The ordering on $\E$ is complete.
\end{Theorem}
\Proof
Let $S$ be a non-empty subset of $\E$ bounded above by $X$, say. I.e., for every $x \in S$,
$x < X$.
We want to exhibit a supremum for $S$.
If  $S$ contains a greatest element, i.e., if there is an $s \in S$ such that $x \le s$ for all $x \in S$,
then that greatest element is the supremum we need (the second part of the definition of the supremum
being vacuous if $s \in S$).
So we may assume that for any $x \in S$ there is a $y \in S$ with $x < y$.

Define a function $f$ from $\Z$ to $\Z$ as follows:
\[
f(p) =
 \left\{
    \begin{array}{ll}
	\Max{\{\Floor{\Eu{p}x}|x \in S\}} & \mbox{if $p \ge 0$} \\
	-f(-p)	& \mbox{otherwise}
    \end{array}
 \right.
\]
This is well-defined, since the sets of integers whose maxima we take are
certainly non-empty, since $S$ is non-empty, and these sets are bounded above, since if $x \in S$,
we have $x < X$ whence $\Floor{px} < \Floor{pX} + 1$ for $p \ge 0$.

I claim $f$ is an almost homomorphism. By lemma~\ref{ExtensionLemma}, it suffices to show
that $d_f(m, n) = f(m+n) - f(m) - f(n)$ is bounded as $m$ and $n$ range over $\N$.
By the definition of $f$, for any $m \in \N$, there is an $x_m\in S$ such that $f(m) = \Floor{\Eu{m}x_m}$
and for any $y \in S$, $\Floor{\Eu{m}y} \le \Floor{\Eu{m}x_m}$. Given $m$ and $n$ in $\N$, let $x$ be
the largest of $x_m$, $x_n$ and $x_{m+n}$, so that $f(m) = \Floor{\Eu{m}x}$,
$f(n) = \Floor{\Eu{n}x}$ and $f(m+n) = \Floor{\Eu{(m+n)}x} = \Floor{\Eu{m}x + \Eu{n}x}$.
By lemma~\ref{FloorAHLemma}, it follows that $d_f(m , n) = \Floor{\Eu{m}x + \Eu{n}x} - \Floor{\Eu{m}x} - \Floor{\Eu{n}x}$ is bounded.

I claim $s = [f]$ is the supremum of $S$.
To prove this we must show, {\em(i)}, that $x \le s$ whenever $x \in S$, and, {\em(ii)}, that if $y$ in $\E$
is such that $x \le y$ whenever $x \in S$ then $s \le y$.

For {\em(i)}, let $x \in S$. Since we are assuming that $S$ has no greatest element, there is a $y \in S$
with $x < y$. By lemma~\ref{ArchLemma}, there are integers $M$ and $N$ with $N$ positive such that
$x < \Eu{M}/\Eu{N} < y$, so that $\Eu{N}x < \Eu{M} < \Eu{N}y$.
Now for any positive $p$, we have the following (in $\Z$):
\[
f(pN) = \Max{\{\Floor{\Eu{(pN)}x}|x \in S\}} \ge \Floor{\Eu{(pN)}(\Eu{M}/\Eu{N})} = \Floor{\Eu{(pM)}} = pM
\]

Where the inequality holds because $y \in S$ and $\Floor{\Eu{(pN)}y} \ge \Floor{\Eu{(pN)}(\Eu{M}/\Eu{N})}$ by our choice
of $M$ and $N$.
It follows from the definitions of the ordering and multiplication that (in $\E$):
\[
s\Eu{N} = [f]\Eu{N} = [p \mapsto f(pN)] \ge [p \mapsto pM] = \Eu{M}
\]
But then we have $s \ge \Eu{M}/\Eu{N} > x$ completing the proof of {\em(i)}.

For {\em(ii)}, let $y \in \E$ be such that $x \le y$ whenever $x \in S$.
Assume for a contradiction that $y < s$. As in the proof of {\em(i)}, there are then integers $M$ and $N$ with $N$ positive such that  $\Eu{N}y < \Eu{M} < \Eu{N}s$ and so, as $\Eu{N}s = s\Eu{N} = [p \mapsto f(pN)]$
and $\Eu{M} = [p \mapsto pM]$, from the definition of the ordering we have that $f(pN) > pM$
for infinitely many positive $p$. But this is impossible, since for any $p$, there is an $x$ in $S$
such that $f(pN) = \Floor{\Eu{(pN)}x}$, but $\Eu{(pN)}x < \Eu{(pN)}y < \Eu{pM}$ and 
so $f(pN) = \Floor{\Eu{(pN)}x} \le \Floor{\Eu{(pM)}} = pM$.
This completes the proof of {\em(ii)} and hence the proof of the theorem.
\Done

Theorems~\ref{OrderedFieldThm} and~\ref{CompleteThm} together show that $\E$ is
a complete ordered field, and, as we have already remarked, this means that $\E$
is isomorphic to the the ordered field of real numbers $\R$, as constructed, for example,
using the method of Dedekind cuts, or the Cantor-Meray method using Cauchy sequences. 

\section{Sources and Remarks}

The theorem that any two complete ordered fields are isomorphic is essentially
due to H\"{o}lder \cite{Hoelder01}. A modern account may be found in~\cite{Ebbinghaus90}.
In fact, any two non-trivial dense complete ordered groups are isomorphic and so it is not
necessary to define the multiplicative structure using representatives: it can be defined by a general consideration
of the order-preserving homomorphisms of the additive structure (e.g., see~\cite{Behrend56}). However, this approach is unattractive for the Eudoxus reals, since the multiplicative structure is very helpful in proving
density and completeness.

The Eudoxus real number representation is due to Stephen Schanuel and dates back to the early 1980s.
Schanuel observed that the graph of the function $p \mapsto \Floor{xp}$ is a subset of $\Z \times \Z$
that can be thought of a discrete representation of the real quantity $x$. He named the
resulting development of the real numbers after Eudoxus, since it seemed to reflect the relationship
between the discrete and the continuous apparent in the ancient theory of proportion.

Schanuel communicated the idea to many people, but did not publish it.
The Eudoxus reals have cropped up from time to time
in various mathematical forums on the Internet over the last decade or so,
but as far as I am aware, the only descriptions in the literature are:

\begin{itemize}
\item
a short article by Ross Street~\cite{Street85};
\item
chapter 2 of John Harrison's thesis \cite{Harrison96};
\item
a recent report of an independent discovery by Norbert A'Campo \cite{ACampo03}.
\end{itemize}

All of these references use the ideas to construct the reals from first principles.
Street cites Schanuel as his source for the construction and mentions an independent discovery
by Richard Lewis.
Harrison learnt of the ideas of Schanuel and developed them into a construction of the reals suitable
for a mechanized implementation of the basic of analysis using the HOL theorem prover.
A'Campo remarks that his work has its origins on Poincar\'{e}'s definition of the rotation number of an
orientation-preserving homeomorphism of the circle and gives many other interesting references.

Unfortunately, Street's article leaves many details to the reader and the proof of completeness is not
correct. Street says that the infimum of a non-empty set $S$ of positive reals is represented by $f$ where $f(m) = \Min\{g(m) | [g] \in S\}$, but this will not do: to see this, represent
the set $S$ of positive integers by the almost homomorphisms $g_m, m = 1, 2 \ldots,$
where $g_m(p) = 0$ if $\Abs{p} \le m$, $g_m(p) = mp$ if $\Abs{p} > m$; if $f(m) = \Min\{g_m(m) | m = 1, 2, \ldots\}$,
$[f] = 0$, but the infimum of $S$ is $1$.
Street became aware of the problem shortly after publication of the article and
saw that it could be solved using normal forms for the Eudoxus real numbers.
An improved treatment including a development of
a normal form has been given by his students Ben Odgers and Nguyen Hanh Vo \cite{Odgers02}.

Harrison's treatment deals with a variant of the construction
in which the natural numbers are used to construct the positive real numbers and then the remaining
real numbers are constructed from those. This is very appropriate given the context in which Harrison
was developing the theory, but makes some of the work a little harder than it need be.

I began to write the present note, unaware of the accounts by A'Campo and Street,
simply with a desire to understand the construction
and to provide a sketch of the generalisations given in the appendix.
I had originally expected to follow Harrison fairly closely, but ended
up taking a different route through the maze.
Harrison's approach is based heavily on the observation that an almost homomorphism can
be viewed as a representation of a Cauchy sequence with a linear bound on its rate
of convergence. I found it more natural to focus on the almost homomorphism property itself.
The result turns out to be more like a rehabilitation of the sketch given by Street.

The problem with Street's original definition of infima is an instance of a general feature of
the Eudoxus representation: no finite amount of information about an almost
homomorphism $f$ reveals any information about $[f]$, and {\em vice versa}, i.e., for any $x$ in $\E$,
any finite partial function from $\Z$ to $\Z$ can be extended to an almost homomorphism
that represents $x$.  There are two potential solutions to this: either, as in the present paper and
in Harrison's approach, develop estimates on the growth rates of almost homomorphisms allowing
information about $[f]$ to be obtained from global information about $f$ or, as in Street's
proposal, to normalise $f$ so that finite information about it does lead to information about
$[f]$. A'Campo gives an ingenious and simple construction of a very efficient normal form
for which $\Abs{d_f(x, y)} \le 1$.  

Finally, it is worth pointing out the normal form for the Eudoxus reals that arises
naturally from the approach of the present note: if $0 < x \in  \E$, we can define
a function $\nu_x$ from $\Z$ to $\Z$ by $\nu_x(p) = \Floor{\Eu{p}x}$. 
By our proof of completeness, one has $x = \Sup{\{x\}} = [\nu_x]$.
The reader who cares to calculate
the first few values of $\nu_{\Rii}$ may well discern an uncanny resemblance between
the results and figure~\ref{Colonnade}.

\section*{Acknowledgments}
I am grateful to Ross Street and Stephen Schanuel for information about the
history of the Eudoxus reals and to Norbert A'Campo for help with
my understanding of his formulation.

\bibliographystyle{plain}
{\raggedright
\bibliography{bookspapers}
}

\vfill
{\tiny RCS Info: $Revision: 1.25 $ $Date: 2004/05/24 15:02:14 $}

\newpage
\appendix
\section{Generalisations}
Joe Shipman has observed that the group $\Z$ in the definition of an almost homomorphism
can be replaced by an arbitrary abelian group $G$, if we consider a subset of $G$ to be bounded
iff. it is finite. In fact, it is convenient to generalise
a little further: given abelian groups $G$, $H$, the set of all functions $G \rightarrow H$, is an
abelian group under pointwise addition of functions, and there is a chain of subgroups:
\[
B(G, H) \prec S(G, H) \prec (G \rightarrow H)
\]
where $B(G, H)$ is the set of functions with finite range and $S(G, H)$ is the set of almost
homomorphisms, i.e., functions $f$ for which $d_f(x, y) = f(x + y) - f(x) - f(y)$ defines a function with
finite range.
So for example, by lemma~\ref{FloorAHLemma} above, the floor function is an element
of $S(\R, \Z)$.
Let us write $\E(G, H)$ for the quotient $S(G, H)/B(G, H)$, so that
the group $\E$ of definition~\ref{Edef} is $\E(\Z, \Z)$.

If one of the composable functions $f$ and $g$ is a genuine homomorphism and
the other is an almost homomorphism, then the composite $f \circ g$ is also an almost homomorphism.
Writing $\Ab$ for the category of abelian groups,
it follows that $\E(\_, \_)$ is a bifunctor $\Ab \times \Ab \rightarrow \Ab$
contravariant in its first argument and covariant in its second.
One may then check%
\footnote{The only non-trivial detail to take care of is in the left-to-right direction of the 
first equivalence, where one needs to show  that $(x, y) \mapsto  f(x, 0) + f(0, y)$
and  $f$ represent the same thing in $\E(G_1 + G_2, H)$, but that is an 
immediate consequence of $f$ being an almost homomorphism.}
that $\E(\_, \_)$ commutes with finite direct sums and products in the sense that
there are natural equivalences:
\[
    \begin{array}{rcl}
\E(G_1 + G_2, H) &\cong& \E(G_1, H) + \E(G_2, H)\\
\E(G, H_1 \times H_2) &\cong& \E(G, H_1) \times \E(G, H_2)
    \end{array}
\]

If either of $G$ or $H$ is finite, then $\E(G, H) = 0$. Since we have shown that $\E(\Z, \Z) \cong \R$,
it follows from the basis theorem for finitely generated abelian groups
that, if $G$ and $H$ are finitely generated, then there is a natural equivalence
\[
    \begin{array}{rcl}
\E(G, H) &\cong& \HomAb{G \otimes \R}{H \otimes \R}
    \end{array}
\]
which confirms Shipman's conjecture that $\E(\Z^n, \Z^n)$ is isomorphic
to the group of $n \times n$ matrices over $\R$ under addition.

By further adapting our notion of boundedness, other possibilities arise. For an example,
consider continuous functions between topological abelian groups where a function
is taken to be bounded iff. the closure of its range is compact. The definitions above then correspond to
the special case where the topology is discrete. Under this generalisation, using the usual
topology on $\R$, we have $\E(\R, \R) \cong \E(\Z, \Z) \cong \R$, where the  isomorphism $\E(\Z, \Z) \cong \E(\R, \R)$
is induced by mapping an almost homomorphism in $\Z \rightarrow \Z$ to its piecewise linear extension to an almost
homomorphism in $\R \rightarrow \R$.

\end{document}